\documentclass[a4,12pt]{article}

\usepackage[english]{babel}

\usepackage{amsfonts}
\usepackage{amsmath}
\usepackage{amssymb}
\usepackage{graphicx}

\usepackage{array}
\usepackage{color}

\newcommand{\Argmin}{\mathop{\mathrm{Arg}\!\min}}

 \newcommand{\RR}{{\mathbb R}}
\newtheorem{remark}{Remark}
\newtheorem{theorem}{Theorem}
\newtheorem{Coro}{Corollary}

\title{
\begin{center}
{Workshop ``Optimization and Statistical Learning''\\ April 10--14, 2017, Les Houches, France
}
\end{center}
\bigskip
\bigskip
\textbf{Algorithms of Inertial Mirror Descent\\ in Convex Problems of Stochastic Optimization%
\footnote{The full paper is accepted at Russian journal \textit{Automatika i Telemekhanika}
which would be translated as \textit{Automation and Remote Control}.}
  }
}

\date{April 12, 2017}

\author{{}
\mbox{}\vspace{7pt}
\textsc{Alexander Nazin}
        }

\begin{document}

\maketitle
\begin{center}
\vspace{-17pt}

\textsf{ICS RAS, Moscow, Russia
}
\end{center}
\vspace{7pt}

\abstract{
~The goal is to modify the known method of mirror descent (MD), proposed by A.S.~Nemirovsky and D.B.~Yudin in 1979.
The paper shows the idea of a new, so-called inertial MD method with the example of a deterministic optimization problem in continuous time.
In particular, in the Euclidean case, the heavy ball method by B.T.~Polyak is realized.
It is noted that the new method does not use additional averaging.
A discrete algorithm of inertial MD is described. The theorem on the upper bound on the error in the objective function is proved.%
}
\vspace{7pt}
\\

\centerline{\textbf{1. The idea of method of inertial mirror descent
            }}
Let $f:\mathbb{R^N}\to\mathbb{R}$ be convex, differentiable function having a unique minimum point
$x^{*}\in\mathrm{Argmin}f(x)$ and its minimal value
 $f^{*}=f(x^{*})$.
 Consider 
 continuous algorithm which extends MDM that is
\begin{eqnarray}
  \dot{\zeta}(t) &=& - \nabla f(x(t)),\quad t\geq0,\,\,\zeta(0)=0,
  \label{dotzeta}
  \\
  \mu_t\dot{x}(t) + x(t) &=& \nabla W(\zeta(t)),\quad x(0)=\nabla W(\zeta(0))
    \label{mudotx}.
\end{eqnarray}
Functional parameter in \eqref{mudotx} is a convex, continuously differentiable function
$W:\mathbb{R^N}\to\mathbb{R}_+$ having conjugate function
\begin{equation}\label{Legendre}
V(x)=\sup_{\zeta\in\mathbb{R^N}}{\{\langle\zeta,x\rangle - W(\zeta)\}}.
\end{equation}
Let $W(0)=0$
, %
$V(0)=0$, and
 $\nabla W(0)=0$ for simplicity.

\begin{remark}
Under parameter
 $\mu_t\equiv 0$ in \eqref{mudotx}, algorithm \eqref{dotzeta}--\eqref{mudotx}
represents MDM (in continuous time) \cite{NemirovskiEtAl-1979};
in particular, the identical map
$\nabla W(\zeta)\equiv\zeta$ and $\mu_t\equiv 0$ lead to a standard gradient method
\begin{equation*}
  \dot{x}(t) = - \nabla f(x(t)),\quad t\geq0.
\end{equation*}
Under $\mu_t\equiv\mu>0$ and $W(\zeta)\equiv\zeta$, algorithm \eqref{dotzeta}--\eqref{mudotx}
leads to continuous method of heavy ball (MHB)
\cite{Polyak-book1983}
\begin{equation*}
  \mu\ddot{x}(t) + \dot{x}(t) = - \nabla f(x(t)),\quad t\geq0.
\end{equation*}
\mbox{}\hfill$\square$
\end{remark}

Further, we assume that differentiable parameter $\mu_t\geq0$, and method
\eqref{dotzeta}--\eqref{mudotx} we call Method of Inertial Mirror Descent (MIDM).

Assume a solution $\{x(t)$, $t\geq0\}$ to system equations \eqref{dotzeta}--\eqref{mudotx} exists.
\\
Consider function
\begin{equation}\label{Lyap}
  W_*(\zeta) =  W(\zeta) - \langle\zeta, x^{*}\rangle\,,\quad   \zeta\in\mathbb{R^N}\,,
\end{equation}
attempting to find a candidate Lyapunov function.

Trajectory derivative to system
\eqref{dotzeta}--\eqref{mudotx}
be
\begin{eqnarray}
  \frac{d}{dt}W_*(\zeta(t)) &=& \langle\dot{\zeta},\nabla{W}  - x^{*} \rangle 
   = -\langle \nabla f(x), \mu_t\dot{x}+x - x^{*} \rangle \leq \label{CandLyap2} \\
   &\leq& f(x^{*})-f(x(t))  -\mu_t\frac{d}{dt}[f(x(t))-f^*]\label{CandLyap3}
\end{eqnarray}
where last inequality results from convexity $f(\cdot)$.
Now, integrating on interval $[0,t]$
with $W_*(0)=0$, we obtain
  \begin{equation}\label{Integrationt}
  \int_0^t f(x(t))dt\, -\, tf^* \leq  -W_*(\zeta(t)) - \mu_t [f(x(t)) -f^*]\Big{|}_0^t + \int_0^t [f(x(t)) -f^*]\dot{\mu}_t dt,
\end{equation}
where two last terms in RHS got by integrating in parts. Taking
\eqref{Legendre} into account, we continue \eqref{Integrationt}:
 \begin{eqnarray}
  \int_0^t f(x(t))dt -tf^* &\leq&  V(x^*) - \mu_t [f(x(t)) -f^*] +\nonumber \\
  && +\,\, \mu_0 [f(x(0)) -f^*] +\nonumber\\
  && +\, \left[\sup_{s\in[0,t]}{\dot{\mu}_s}\right] \int_0^t [f(x(t)) -f^*] dt. \label{IntegrationtCont2}
\end{eqnarray}
Therefore, it is reasonable to introduce the following constraints on patameter
$\mu_t\geq0$ :
\begin{equation}\label{mu01}
  \mu_0=0,\quad \dot{\mu}_t\leq1\,\,\forall t>0,
\end{equation}
leading to inequality
\begin{equation*}
  f(x(t)) -f^* \leq  V(x^*)/\mu_t\,.
\end{equation*}
Maximizing $\mu_t$ under constraints \eqref{mu01} we get
\begin{equation}\label{muteqt}
  \mu_t = t, \quad t\geq0.
\end{equation}
The related (continuous) IMD algorithm
\begin{eqnarray}
  \dot{\zeta}(t) &=& - \nabla f(x(t)),\quad t\geq0,\quad \zeta(0)=0,
  \label{dotzetat}
  \\
  t\,\dot{x}(t) + x(t) &=& \nabla W(\zeta(t)),
    \label{mudotxt}
\end{eqnarray}
proves upper bound
\begin{equation}\label{simplerUBt2}
  f(x(t)) -f^* \leq 
                     V(x^*)\,t^{-1},\quad\forall t>0\,.
\end{equation}
\\

\centerline{\textbf{2. Stochastic optimization problem
            }}
Consider minimization problem
\begin{equation}\label{Afunc}
f(x) \triangleq \mathbb{E}\, Q(x,Z)\,\to
\min_{x\in X}\,,
\end{equation}
where {\rm loss function}
$Q:X\times\mathcal{Z}\to\mathbb{R}_+$\,
contains random variable $Z$ with unknown
distribution on space $\mathcal{Z}$,
$\mathbb{E}$ --- mathematical expectation,
set $X\subset\RR^N$ --- given convex compact in $N$-dimension space,
random function\, $Q(\cdot\,,Z):X\to\mathbb{R}_+$\, is convex a.s. on $X$.

Let i.i.d sample $(Z_1,\dots,Z_{t-1})$ be given
where all $Z_i$ have the same distribution on $\mathcal{Z}$ as $Z$.
%
Introduce notation for stochastic subgradients
\begin{equation}\label{Genu}
u_k(x) = \nabla_{x}^{}Q(x,Z_k )\,, \quad k=1,2,\dots,
\end{equation}
such that
$\forall x\in X$,
$$\mathbb{E}\,u_k(x) \in \partial f(x).$$
The goal is in constructing and proving novel recursive MD algorithms
meant for minimization \eqref{Afunc}
 and using stochastic subgradients $u_t(x_{t-1})$\, \eqref{Genu} at current points $x=x_{t-1}\in X$, $t=1,2,\dots$.
\vspace{12pt}

\centerline{\textbf{3. Algorithm IMD.
Main results.}}

Let $\|\cdot\|$ be a norm in primal space $E=\mathbb{R}^N$, and $\|\cdot\|_*$
be the related norm in dual space $E^*=\mathbb{R}^N$;
set $X\subset E$ is convex compact.

{\bf Assumption
(L).} \textit{Convex function $V: {X}\to \mathbb{R}_+$ is such that its
$\beta$-conjugate  $W_\beta$ is continuously differential on $E^*$ with gradient $\nabla W_{\beta}$
satisfying Lipschitz condition
\[
\| \nabla W_{\beta}(\zeta) - \nabla W_\beta(\,\tilde{\zeta}\,) \|
 \le
\frac{1}{\alpha\beta} \|\zeta-\tilde{\zeta} \|_*^{} \,,
\quad
\forall\, \zeta, \tilde{\zeta} \in E^*,  \;\beta>0,
\]}
\textit{where $\alpha$ is positive constant being independent of $\beta$.}

Consider now the discrete time $t\in Z_+$.
Write a discrete version of algorithm IMD \eqref{dotzetat}--\eqref{mudotxt}
using stochastic subgradients \eqref{Genu} instead of the gradients $\nabla f(\cdot)$:
\begin{eqnarray}
  \tau_t &=& \tau_{t-1} + \gamma_t\,,\quad t\geq1,\, \tau_0=0, \label{algtk} \\
  \zeta_t &=& \zeta_{t-1} + \gamma_t u_t(x_{t-1}),\quad \zeta_{0}=0, \label{algzk} \\
  \tau_t\,\frac{x_t - x_{t-1}}{\gamma_{t+1}} + x_{t} &=& -\nabla W_{\beta_{t}}(\zeta_t),\quad x_0=-\nabla W_{\beta_{0}}(\zeta_0).
  \label{algxk}
\end{eqnarray}
Here function $W_\beta$ is defined by proxy-function $V:X\to\mathbb{R}_+$
via Legendre--Fenchel transformation, i.e.
\begin{equation}\label{fench}
W_\beta(\zeta) = \sup_{x\in X} \left\{ -\zeta^T x -\beta V(x)\right\}%
, \quad \zeta\in E^*
\, ,
\end{equation}

\begin{remark}
Equation \eqref{algxk} may be written as
\begin{equation}\label{algxkSled}
      x_t = \frac{\tau_t}{\tau_t +\gamma_{t+1}} x_{t-1} - \frac{\gamma_{t+1}}{\tau_t +\gamma_{t+1}} \nabla W_{\beta_{t}}(\zeta_t).
\end{equation}
Since the vectors\, $[-\nabla W_{\beta_{t}}(\zeta_t)]\in X$ under each $t\geq0$, equations
\eqref{algtk}--\eqref{algzk} show that $x_t\in X$ by induction.
 \hfill$\Box$
\end{remark}

Further, let sequences $(\gamma_i)_{i\ge
1}$ and $(\beta_i)_{i \ge 1}$ are of view
\begin{equation}\label{gb1}
\gamma_i\equiv 1\,,\quad \beta_i =\beta_0\sqrt{i+1}\,,\quad i=1,2,\dots,
\quad \beta_0>0.
\end{equation}
Then system equations \eqref{algtk}--\eqref{algxk} leads to the IMD algorithm:
\begin{eqnarray}
  \zeta_t &=& \zeta_{t-1} +  u_t(x_{t-1}),\quad \zeta_{0}=0,\quad x_0=-\nabla W_{\beta_{0}}(\zeta_0), \label{algzkt} \\
    x_t &=&  x_{t-1} - \frac{1}{t+1}\left( x_{t-1} + \nabla W_{\beta_{t}}(\zeta_t)\right),\quad t\geq1
    .
  \label{algxt}
\end{eqnarray}

\begin{theorem}\label{th:1}
Let ${X}$ be convex closed set in $\mathbb{R}^N$,
and loss function $Q(\cdot,\cdot)$ satisfies the conditions of section~2,
and, moreover,
\begin{equation}
\label{eq:supth2}
\sup_{{x}\in{X}}\mathbb{E}\|\nabla_{{x}}^{}\,Q({x},Z)\|_{*}^{2}\leq
L_{{X},\,Q}^2
\,,
\end{equation}
where constant $L_{{X},\,Q}\in(0,\infty)$.
Let $V$ be proxy-function on ${X}$ with parameter $\alpha>0$ from assumption
\textrm{(L)},
and let exists minimum point ${x}^* \in\displaystyle
\Argmin_{{x}\in{X}} f({x})$. Then for any $t\ge 1$ estimate ${x}_t$\,, defined by algorithm
\eqref{algzkt}, \eqref{algxt} with stochastic subgradients (\ref{Genu})
and sequence $(\beta_i)_{i\ge 1}$ from
(\ref{gb1}) with arbitrary $\beta_0 >0$, satisfies inequality
\begin{eqnarray}
\nonumber
\mathbb{E}\, f({x}_t) -
\min_{{x}\in{X}}f({x})
 \le
 \left(\beta_0  V({x}^*)+{L_{{X},\,Q}^2
 \over \alpha\beta_0}\right){\sqrt{t+2}\over t+1}\,.
\end{eqnarray}
\hfill$\Box$
\end{theorem}

\begin{Coro}
If constant $\overline V$ in Theorem 1 assumptions is such that $V({x}^*)\le
\overline  V$ and $\beta_0 = L_{{X},\,Q}\,(\alpha\, \overline  V\,)^{-1/2}$
then
\begin{eqnarray}\label{eq:th31}
\mathbb{E}\, f({x}_t) -
\min_{{x}\in{X}}f({x})
 \le
 2\, L_{{X},\,Q}
 \left(\alpha^{-1}\overline  V\,\right)^{1/2}
 {\sqrt{t+2}\over t+1} \, .
\end{eqnarray}
In particular,
one may get
$\overline  V =\displaystyle\max_{{x}\in{X}}
V({x})$.
\hfill$\Box$
\end{Coro}

\newpage
\small %

\end{document}